\documentclass[]{article}

\usepackage[utf8]{inputenc}

\usepackage{amsmath}
\usepackage{amssymb}
\usepackage{latexsym}
\usepackage{bbm}
\usepackage{graphicx}
\usepackage{url}

\newcommand{\R}{\mathbb{R}}
\newcommand{\N}{\mathbb{N}}
\DeclareMathOperator{\convOp}{conv}
\newcommand{\conv}[1]{\convOp({#1})}
\DeclareMathOperator{\cconeOp}{ccone}
\newcommand{\ccone}[1]{\cconeOp({#1})}
\DeclareMathOperator{\affOp}{aff}
\newcommand{\aff}[1]{\affOp({#1})}
\DeclareMathOperator{\matchOp}{\mathcal{M}}
\newcommand{\allMatch}[1]{\matchOp({#1})}
\newcommand{\match}[2]{\matchOp_{#1}({#2})}
\newcommand{\charVec}[1]{\chi({#1})}
\DeclareMathOperator{\polyOp}{P}
\newcommand{\allMatchPoly}[1]{\polyOp^{\text{match}}({#1})}
\newcommand{\matchPoly}[2]{\polyOp^{\text{match}}_{{#1}}({#2})}
\newcommand{\sptPoly}[1]{\polyOp^{\text{spt}}({#1})}
\newcommand{\knapPoly}[2]{\polyOp^{\text{knap}}({#1},{#2})}
\newcommand{\setDef}[2]{\{{#1}\,:\,{#2}\}}
\newcommand{\outEdges}[1]{\delta({#1})}
\newcommand{\perm}[1]{\polyOp^{\text{perm}}({#1})}
\newcommand{\ints}[1]{[{#1}]}
\newcommand{\birk}[1]{\polyOp^{\text{birk}}({#1})}

\newcommand{\scalProd}[2]{\langle{#1},{#2}\rangle}
\newcommand{\facLat}[1]{\mathcal{L}({#1})}
\newcommand{\unitVec}[1]{\mathbbm{e}_{#1}}
\newcommand{\zeroVec}[1]{\mathbb{O}_{#1}}
\newcommand{\RNonNeg}{\R_+}
\DeclareMathOperator{\extForm}{EF}
\DeclareMathOperator{\bigOOp}{O}
\newcommand{\bigO}[1]{\bigOOp({#1})}

\newtheorem{theorem}{Theorem}

\title{Extended Formulations in Combinatorial Optimization}
\author{Volker Kaibel\thanks{Institut f\"ur Mathematische Optimierung, Fakult\"at f\"ur Mathematik, Otto-von-Guericke Universit\"at Magdeburg, Universit\"atsplatz~2, 39108~Magdeburg, Germany, \texttt{kaibel@ovgu.de}}}

\date{\today}

\begin{document}

\maketitle

\begin{abstract}
	The concept of representing a polytope that is associated with some combinatorial optimization problem as a linear projection of a higher-dimensional polyhedron has recently received increasing attention. In this paper (written for the newsletter \emph{Optima} of the Mathematical Optimization Society), we provide a brief introduction to this topic  and sketch some of the recent developments with respect to both tools for constructing such  extended formulations as well as  lower bounds on their sizes. 
\end{abstract}


\section{Introduction}
Linear Programming based methods and polyhedral theory form the backbone of large parts of Combinatorial Optimization. The basic paradigm here is to identify the feasible solutions to a given problem with some vectors in such a way that the optimization problem becomes the problem of optimizing a linear function over the finite set~$X$ of these vectors. The optimal value of a linear function over~$X$ is equal to its  optimal value over the convex hull~$\conv{X}=\setDef{\sum_{x\in X}\lambda_xx}{\sum_{x\in X}\lambda_x=1,\lambda\ge\zeroVec{}}$ of~$X$. 
According to the Weyl--Minkowski Theorem~\cite{Wey35,Min96}, every \emph{polytope} (i.e., the convex hull of a finite set of vectors) 
can be written as the set of solutions to a system of linear equations and inequalities. Thus one  ends up with a linear programming problem. 

As for the maybe most classical example, let us consider the set $\allMatch{n}$ of all matchings in the complete graph $K_n=(V_n,E_n)$ on~$n$ nodes (where a matching is a subset of edges no two of which share a common end-node). Identifying every matching $M\subseteq E_n$ with its characteristic vector $\charVec{M}\in\{0,1\}^{E_n}$ (where $\charVec{M}_e=1$ if and only if $e\in M$), we obtain the \emph{matching polytope} 
\begin{equation*}
	\allMatchPoly{n}=\convOp\setDef{\charVec{M}}{M\in\allMatch{n}}\,.
\end{equation*}
In one of his seminal papers, Edmonds~\cite{Edm65b} proved that $\allMatchPoly{n}$ equals the set of all $x\in\RNonNeg^{E_n}$ that satisfy the inequalities $x(\outEdges{v})\le 1$ for all $v\in V_n$ and
	$x(E_n(S))\le \lfloor |S|/2\rfloor$ for all subsets $S\subseteq V_n$ of \emph{odd} cardinality $3\le |S|\le n$ (where $\outEdges{v}$ is the set of all edges incident to~$v$, $E_n(S)$ is the set of all edges with both end-nodes in~$S$, and $x(F)=\sum_{e\in F}x_e$). No inequality in this system, whose size is exponential in~$n$, is redundant.

The situation is quite similar for the \emph{permutahedron}~$\perm{n}$, i.e., the convex hull of all vectors that arise from permuting the components of $(1,2,\dots,n)$. Rado~\cite{Rad52} proved
that~$\perm{n}$ is described by the equation $x(\ints{n})=n(n+1)/2$ and the inequalities $x(S)\ge |S|(|S|+1)/2$ for all $\varnothing\ne S\subsetneq\ints{n}$ (with $\ints{n}=\{1,\dots,n\}$), none of the  $2^n-2$ inequalities being redundant. However if for each  permutation $\sigma:\ints{n}\rightarrow\ints{n}$ we consider the corresponding \emph{permutation matrix} $y\in\{0,1\}^{n\times n}$ (satisfying $y_{ij}=1$ if and only if $\sigma(i)=j$) rather than the vector 
$(\sigma(1),\dots,\sigma(n))$, we obtain a much smaller description of the resulting polytope, since,
according to Birkhoff~\cite{Bir46} and von Neumann~\cite{vNeu53}, the convex hull~$\birk{n}$ (the \emph{Birkhoff-Polytope}) of all  $n\times n$-permutation matrices is equal to the set of all \emph{doubly-stochastic} $n\times n$-matrices (i.e., nonnegative $n\times n$-matrices all of whose row- and column sums are equal to one). It is easy to see that the permutahedron~$\perm{n}$
 is a linear projection of the Birkhoff-polytope~$\birk{n}$ via the map defined by $p(y)_i=\sum_{j=1}^n jy_{ij}$. Since, for every linear objective function vector $c\in\R^n$, we have $\max\setDef{\scalProd{c}{x}}{x\in\perm{n}}=\max\setDef{\sum_{i=1}^n\sum_{j=1}^njc_iy_{ij}}{y\in\birk{n}}$,  one can use~$\birk{n}$ (that can be described by the $n^2$ nonnegativity-inequalities) instead of~$\perm{n}$ (whose description requires $2^n-2$ inequalities) with respect to  linear programming related issues. 

In general, an \emph{extension} of a polytope~$P\subseteq\R^n$  is a polyhedron~$Q\subseteq\R^d$ (i.e., an intersection of finitely many affine hyperplanes and halfspaces) together with a linear projection $p:\R^d\rightarrow\R^n$ satisfying $P=p(Q)$. Any description of~$Q$ by linear equations and linear inequalities  then (together with~$p$) is an \emph{extended formulation} of~$P$. The \emph{size} of the extended formulation is the number of inequalities in the description. Note that we  neither account for the number of equations (we can get rid of them by eliminating variables) nor for the number of variables (we can  ensure that there are not more variables than inequalities by projecting~$Q$  to the orthogonal complement of its \emph{lineality space}, where the latter is the space of all directions of lines contained in~$Q$).
If~$T\in\R^{n\times d}$ is the matrix with $p(y)=Ty$, then, for every $c\in\R^n$, we have $\max\setDef{\scalProd{c}{x}}{x\in P}=\max\setDef{\scalProd{T^tc}{y}}{y\in Q}$. 

In the example described above, $\birk{n}$ thus provides an extended formulation of~$\perm{n}$ of size~$n^2$. It is not known whether one can do something similar for the matching polytopes~$\allMatchPoly{n}$ (we will be back to this question in Section~\ref{subsec:symm}). However there
 are many other examples of nice and small extended formulations for polytopes associated with combinatorial optimization problems. The aim of this article (that has appeared in~\cite{Optima85}) is to show a few of them and to shed some light on the geometric, combinatorial and algebraic background of this concept that recently has received  increased attention. The presentation is not meant to be a survey (for this purpose, we refer to  Vanderbeck and Wolsey~\cite{VW10} as well as to Cornu\'{e}jols, Conforti, and Zambelli~\cite{CCZ10}) but rather  an appetizer for investigating alternative possibilities to express combinatorial optimization problems by means of linear programs. 

While we will not be concerned with practical aspects here, extended formulations have also proven to be  useful in computations. You can find more on this in Laurence Wolsey's discussion column in~\cite{Optima85}. Fundamental work with respect to understanding the concept of extended formulations and its limits has been done by Mihalis Yannakakis in his 1991-paper \emph{Expressing Combinatorial Optimization Problems by Linear Programs}~\cite{Yan91} (see Sect.~\ref{subsec:positiveRank} and~\ref{sec:fundLims}). He discusses some of his thoughts on the subject in another discussion column in~\cite{Optima85}.

\section{Some Examples}
\label{sec:ex}

\subsection{Spanning Trees}
\label{subsec:spt}
The \emph{spanning tree polytope}~$\sptPoly{n}$ associated with the complete graph $K_n=(V_n,E_n)$ on~$n$ nodes is the convex hull of all characteristic vectors of spanning trees, i.e., of all subsets of edges that form connected and cycle-free subgraphs. In another seminal paper, Edmonds~\cite{Edm71} proved that~$\sptPoly{n}$ is the set of all $x\in\RNonNeg^{E_n}$ that satisfy the equation $x(E_n)=n-1$ and the inequalities $x(E_n(S))\le |S|-1$ for all $S\subseteq V_n$ with $2\le |S|<n$. Again, none of the exponentially many inequalities is redundant. 

However, by introducing additional variables $z_{v,w,u}$ for all ordered triples $(v,w,u)$ of pairwise different nodes meant to encode whether the edge~$\{v,w\}$ is contained in the tree and~$u$ is in the component of~$w$ when removing~$\{v,w\}$ from the tree,  it turns out that 
the system consisting of the equations 
$x_{\{v,w\}} - z_{v,w,u}-z_{w,v,u} =0$ and
$x_{\{v,w\}} + \sum_{u\in \ints{n}\setminus\{v,w\}}z_{v,u,w} = 1$
(for all pairwise different $v,w,u\in V_n$)
along with the nonnegativity constraints and the equation $x(E_n)=n-1$ provides an extended formulation of~$\sptPoly{n}$ of size~$\bigO{n^3}$ (with orthogonal projection to the space of $x$-variables). 
This formulation  is due to  Martin~\cite{KM91} (see also~\cite{Yan91,CCZ10}). You will find an alternative one in Laurence Wolsey's discussion column below.

\subsection{Disjunctive Programming}
\label{subsec:disPrg}
If $P_i\subseteq\R^n$ is a polytope for each $i\in\ints{q}$, then clearly $P=\conv{P_1\cup\cdots\cup P_q}$ is a polytope as well, but, in general, it is  difficult to derive a description by linear equations and inequalities in~$\R^n$ from such descriptions of the polytopes~$P_i$. However constructing an extended formulation for~$P$ in this situation is very simple. Indeed suppose that each~$P_i$ is described by a system $A^ix\le b^i$ of~$f_i$ linear inequalities (where, in order to simplify notation, we  assume that  equations are written, e.g., as pairs of inequalities). Then the system
$A^iz^i\le \lambda_ib^i$ for all $i\in\ints{q}$, $\sum_{i=1}^q\lambda_i=1$, $\lambda\ge\zeroVec{}$ with variables $z^i\in\R^n$  for all $i\in\ints{q}$ and $\lambda\in\R^q$ is an extended formulation for~$P$ of size $f_1+\cdots+f_q+q$, where the projection is given by $(z^1,\dots,z^q,\lambda)\mapsto z^1+\dots+z^q$. This has been proved first by Balas (see, e.g.,~\cite{Bal85}), even for polyhedra that are not necessarily polytopes (where in this general case~$P$ needs to be defined as the topological closure of the convex hull of the union). 

\subsection{Dynamic Programming}
\label{subsec:dynPrg}
When a combinatorial optimization problem can be solved by a dynamic programming algorithm, one often  can derive an extended formulation for the associated polytope whose size is roughly bounded by the running time of the algorithm. 

A simple example is the 0/1-Knapsack problem, where we are given a nonnegative integral weight vector $w\in\N^n$, a weight bound~$W\in\N$, and a profit vector $c\in\R^n$, and the task is to solve $\max\setDef{\scalProd{c}{x}}{x\in F(w,W)}$ with $F(w,W)=\setDef{x\in\{0,1\}^n}{\scalProd{w}{x}\le W}$. 
A classical dynamic programming algorithm works by setting up an acyclic directed graph with nodes~$s=(0,0)$, $t$, and $(i,\omega)$ for all $i\in\ints{n}$, $\omega\in\{0,1,\dots,W\}$ and arcs from $(i,\omega)$ to $(i',\omega')$ if and only if $i<i'$ and $\omega'=\omega+w_{i'}$, where such an arc would be assigned length~$c_{i'}$, as well as arcs from all nodes to~$t$ (of length zero). Then solving the 0/1-Knapsack problem is equivalent to finding a longest $s$-$t$-path in this acyclic directed network, which can be carried out in linear time in the number~$\alpha$ of arcs. 

The polyhedron~$Q\subseteq\RNonNeg^{\alpha}$ of all $s$-$t$-flows of value one in that network equals the convex hull of all characteristic vectors of $s$-$t$-paths (due to the total unimodularity of the node-arc incidence matrix), thus it is easily seen to be mapped to 
 the associated \emph{Knapsack-polytope} $\knapPoly{w}{W}=\conv{F(w,W)}$ via the projection given by $y\mapsto x$, where $x_i$ is the sum of all components of~$y$ indexed by arcs pointing to nodes of type $(i,\star)$. As~$Q$ is described by nonnegativity constraints, the flow-conservation equations on the nodes different from~$s$ and~$t$ and the equation ensuring an outflow of value one from~$s$,  these constraints provide an extended formulation for $\knapPoly{w}{W}$ of size~$\alpha$. 

However quite often dynamic programming algorithms can only  be formulated as longest-paths problems in acyclic directed \emph{hyper}graphs with hyperarcs of the type~$(S,v)$ (with a subset~$S$ of nodes) whose usage in the path represents the fact that the optimal solution to the partial problem represented by node~$v$ has been constructed from the optimal solutions to the partial problems represented by the set~$S$.  Martin, Rardin, and Campbell~\cite{KRC90} showed that, under the condition that one can assign appropriate \emph{reference sets} to the nodes, also in this more general situation nonnegativity constraints and flow-equations  suffice to describe the convex hull of the characteristic vectors of the hyperpaths. This generalization allows one to derive polynomial size extended formulations for many of the combinatorial optimization  problems that can be solved in polynomial time by dynamic programming algorithms.

\subsection{Others}

A common generalization of the techniques to construct extended formulations by means of disjunctive programming or dynamic programming is provided by  \emph{branched polyhedral systems} (\emph{BPS})~\cite{KL10}. In this framework, one starts from an acyclic directed graph that has associated with each of its non-sink nodes~$v$ a polyhedron in the space indexed by the out-neighbors of~$v$. From these building blocks, one constructs  a polyhedron in the space indexed by all nodes. Under certain conditions one can derive an extended formulation for the constructed polyhedron from extended formulations of the polyhedra associated with the nodes. 

Some very nice extended formulations have recently been given by Faenza, Oriolo, and Stauffer~\cite{FOS10} for stable set polytopes of claw-free graphs. Here the crucial step is to glue together descriptions of stable set polytopes of certain building block graphs by means of \emph{strip compositions}. One of their constructions can be obtained by applying the BPS-framework, though apparently the most interesting one they have  cannot. 

An asymptotically smallest possible extended formulation of size $\bigO{n\log n}$ for the permutahedron~$\perm{n}$  has been found by Goemans~\cite{Goe09}. His construction relies on the existence of \emph{sorting networks} of size $r=\bigO{n\log n}$ (Ajtai, Koml\'os, and Szemer\'edi~\cite{AKS83}), i.e., sequences $(i_1,j_1),\dots,(i_r,j_r)$ for which the algorithm that in each step~$s$ swaps elements~$a_{i_s}$ and~$a_{j_s}$ if and only if $a_{i_s}>a_{j_s}$ sorts every sequence $(a_1,\dots,a_n)\in\R$ into non-decreasing order. The construction principle of Goemans has been generalized to the framework of \emph{reflection  relations}~\cite{KP11}, which, for instance, can be used to obtain small extended formulations for all \emph{$G$-permutahedra of finite reflection groups~$G$} (see, e.g., Humphreys~\cite{Hum90}), including extended formulations of size $\bigO{\log m}$ of regular~$m$-gons, previously constructed by Ben-Tal and Nemirovski~\cite{BN01}. Another application of reflection relations yields extended formulations of size $\bigO{n\log n}$ for \emph{Huffman-polytopes}, i.e., the convex hulls of the leaves-to-root-distances vectors in rooted binary trees with~$n$ labelled leaves. Note that linear descriptions of these polytopes in the original spaces are very large, rather complicated, and unknown (see Nguyen, Nguyen, and Maurras~\cite{NNM10}). 

The list of combinatorial problems for which small (and nice) extended formulations have been found comprises many others, among them perfect matching polytopes of planar graphs (Barahona~\cite{Bar93}), perfectly matchable subgraph polytopes of bipartite graphs (Balas and Pulleyblank~\cite{BP83}), stable-set polytopes of distance claw-free graphs (Pulleyblank and Shepherd~\cite{PS93}),  packing and partitioning orbitopes~\cite{FK09},  subtour-elimination polytopes (Yannakakis~\cite{Yan91} and, for planar graphs, Rivin~\cite{Riv03}, Cheung~\cite{Che03}), and certain mixed-integer programs (see, e.g., Conforti, di Summa, Eisenbrand, and Wolsey~\cite{CDSEW09}).

\section{Combinatorial, Geometric, and Algebraic Background}

\subsection{Face Lattices}
\label{subsec:facLat}
Any intersection of a polyhedron~$P$ with the boundary hyperplane of some affine halfspace containing~$P$ is called a \emph{face} of~$P$. The empty set and~$P$ itself are considered to be (non-proper) faces of~$P$ as well. For instance, the proper faces of a three-dimensional polytope  are its vertices, edges, and the polygons that make up the boundary of~$P$.
Partially ordered by inclusion, the faces of a polyhedron~$P$ form a lattice~$\facLat{P}$, the \emph{face lattice} of~$P$. The proper faces that are maximal with respect to inclusion are the \emph{facets} of~$P$. Equivalently, the facets of~$P$ are those faces whose dimension is one less than the dimension of~$P$. An irredundant linear description of~$P$ has exactly one inequality for each facet of~$P$. 

\begin{figure}[ht]
	\label{fig:embedding}
	\centerline{\includegraphics[height=7cm]{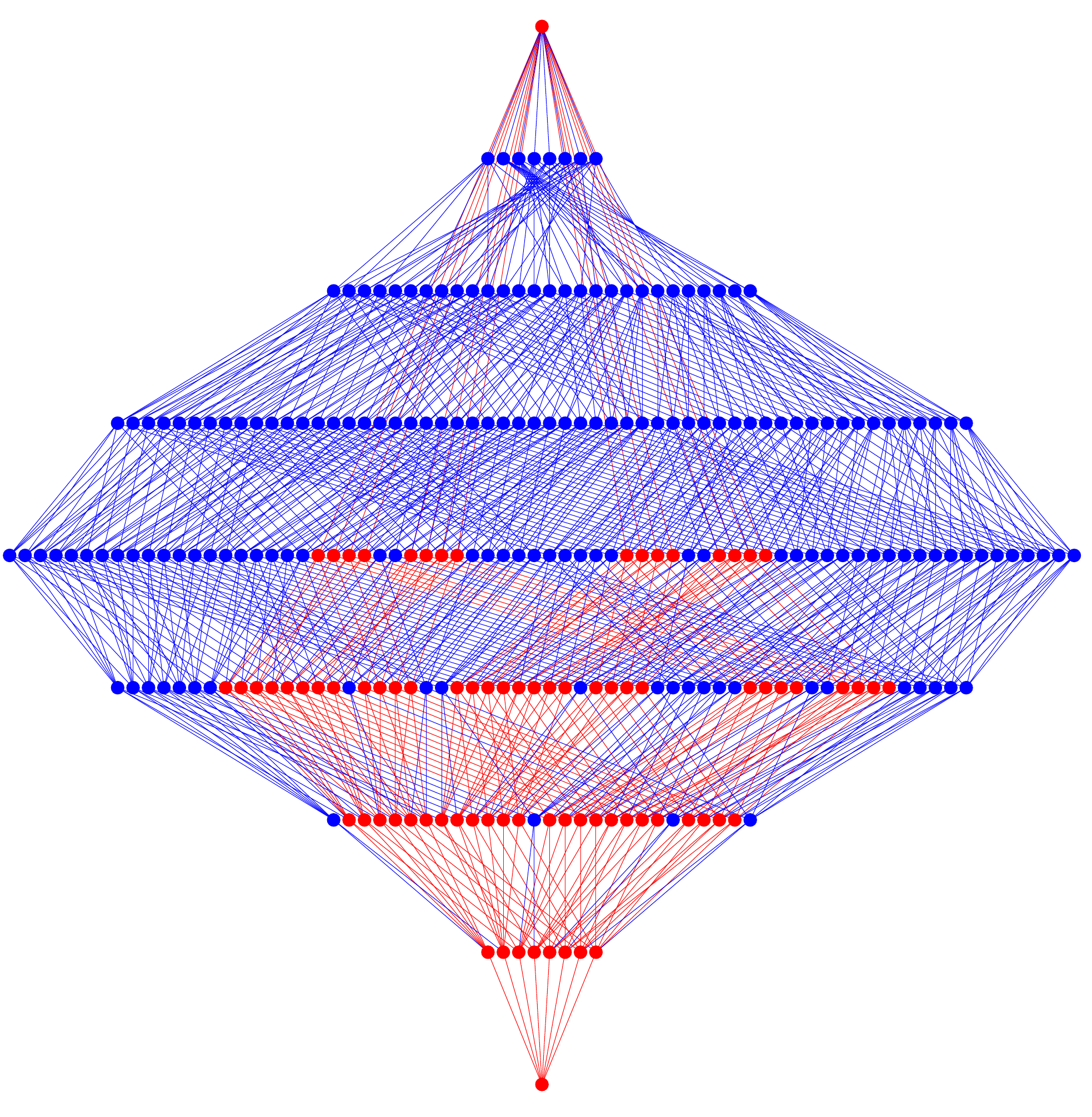}}
	\caption{Embedding of the face lattice of the $4$-dimensional cross-polytope into the face lattice of the $7$-dimensional simplex.}
\end{figure}

If~$Q\subseteq\R^d$ is an extension of the polytope~$P\subseteq\R^n$ with a linear projection~$p:\R^d\rightarrow\R^n$, then mapping each face of~$P$ to its preimage in~$Q$ under~$p$ defines an embedding of~$\facLat{P}$ into~$\facLat{Q}$. Figure~\ref{fig:embedding} illustrates this embedding for the trivial extension $Q=\setDef{y\in\RNonNeg^V}{\sum_{x\in X} y_x=1}$ of $P=\conv{X}$ via $p(y)=\sum_{x\in X}y_xx$ for $X=\{\unitVec{1},-\unitVec{1},\dots,\unitVec{4},-\unitVec{4}\}$ (thus $P$ is the \emph{cross-polytope} in~$\R^4$ with~$16$ facets and~$Q$ is the \emph{standard-simplex} in~$\R^{8}$ with~$8$ facets).
As  this figure suggests, constructing a small extended formulation for a polytope~$P$ means to hide
the  facets of~$P$ in the fat middle part of the face lattice of an extension with few facets.

\subsection{Slack Representations}
Let~$P=\setDef{x\in\mathcal{A}}{Ax\le b}\subseteq\R^n$ be a 
polytope
with affine hull~$\mathcal{A}=\aff{P}$, $A\in\R^{m\times n}$, and $b\in\R^m$. The affine map $\varphi:\mathcal{A}\rightarrow\R^m$ with $\varphi(x)=b-Ax$ (the \emph{slack map} of~$P$ w.r.t. $Ax\le b$) is injective. We denote its inverse (the \emph{inverse slack map}) on its image, the affine subspace  $\tilde{\mathcal{A}}=\varphi(\mathcal{A})\subseteq\R^m$, by $\tilde{\varphi}:\tilde{\mathcal{A}}\rightarrow\mathcal{A}$. The polytope $\tilde{P}=\tilde{\mathcal{A}}\cap\RNonNeg^m$, the \emph{slack-representation} of~$P$ w.r.t. $Ax\le b$, is isomorphic to~$P$ with $\varphi(P)=\tilde{P}$ and $\tilde{\varphi}(\tilde{P})=P$.

If $Z\subseteq\RNonNeg^m$ is a finite set of nonnegative vectors whose \emph{convex conic hull} $\ccone{Z}=\setDef{\sum_{z\in Z}\lambda_zz}{\lambda\ge\zeroVec{}}\subseteq\RNonNeg^m$ contains $\tilde{P}=\tilde{\mathcal{A}}\cap\RNonNeg^m$, then we have $\tilde{P}=\tilde{\mathcal{A}}\cap\ccone{Z}$, and thus, the system $\sum_{z\in Z}\lambda_zz\in\tilde{\mathcal{A}}$ and  $\lambda_z\ge 0$ (for all $z\in Z$) provides an extended formulation of~$P$ of size~$|Z|$ via the projection $\lambda\mapsto\tilde{\varphi}(\sum_{z\in Z}\lambda_zz)$. Let us call such an extension a \emph{slack extension} and the set~$Z$ a \emph{slack generating set} of~$P$ (both w.r.t. $Ax\le b$). 

Now suppose conversely that we have any extended formulation of~$P$ of size~$q$ defining an extension~$Q$ that is \emph{pointed} (i.e., the polyhedron~$Q$ does not contain a line). As  for polytopes above (which in particular are pointed polyhedra), we can consider a 
   slack representation~$\tilde{Q}\subseteq\R^q$ of~$Q$ and the corresponding inverse slack map~$\tilde{\psi}$. Then we have $\varphi(p(\tilde{\psi}(\tilde{Q})))=\tilde{P}$, where~$p$ is the projection map of the extension. If the system $Ax\le b$ is \emph{binding} for~$P$, i.e., each of its inequalities is satisfied at equation by some point from~$P$, then one can show (by using strong LP-duality)
 that there is a \emph{nonnegative} matrix $T\in\RNonNeg^{m\times q}$ with $\varphi(p(\tilde{\psi}(\tilde{z})))=T\tilde{z}$ for all $\tilde{z}\in\tilde{Q}$, thus $\tilde{P}=T\tilde{Q}$. Hence the columns of~$T$ form a slack generating set of~$P$ (w.r.t. $Ax\le b$), yielding a slack extension of size~$q$. As every non-pointed extension of a polytope can be turned into a pointed one of the same size by   projection to the orthogonal complement of the lineality space, we obtain the following result, where the \emph{extension complexity} of a polytope~$P$ is the smallest size of any extended formulation of~$P$.
\begin{theorem}[\cite{FKPT11}]
	\label{thm:extCmplSlackExt}
	The extension complexity of a polytope~$P$ is equal to the minimum size of all slack extensions of~$P$.
\end{theorem}
As 
every slack extension of a polytope is bounded (and since all bounded polyhedra are polytopes), Theorem~\ref{thm:extCmplSlackExt} implies that the extension complexity of a polytope is attained by an extension that is a polytope itself.
Furthermore, in Theorem~\ref{thm:extCmplSlackExt} one may take the minimum over the slack extensions w.r.t. any fixed binding system of inequalities describing~$P$. In particular, all these minima concide.  

\subsection{Nonnegative Rank}
\label{subsec:positiveRank}
Now let~$P=\conv{X}=\setDef{x\in\aff{P}}{Ax\le b}\subseteq\R^n$ be a polytope with some finite set~$X\subseteq\R^n$ and 
$A\in\R^{m\times n}$, $b\in\R^m$. The \emph{slack matrix} of~$P$ w.r.t.~$X$ and~$Ax\le b$ is $\Phi\in\RNonNeg^{\ints{m}\times X}$ with $\Phi_{i,x}=b-\scalProd{A_{i,\star}}{x}$. Thus the slack representation~$\tilde{P}\subseteq\R^m$ of~$P$ (w.r.t. $Ax\le b$) is the convex hull of the columns of~$\Phi$. Consequently, if the columns of a nonnegative matrix $T\subseteq\RNonNeg^{\ints{m}\times \ints{f}}$ form a slack generating set of~$P$, then there is a nonnegative matrix $S\in\RNonNeg^{\ints{f}\times X}$ with $\Phi=TS$. Conversely, for every  factorization $\Phi=T'S'$ of the slack matrix into nonnegative matrices $T'\in\RNonNeg^{\ints{m}\times \ints{f'}}$ and $S'\in\RNonNeg^{\ints{f'}\times X}$, the columns of~$T'$  form a slack generating set for~$P$. 

Therefore constructing an extended formulation of size~$f$ for~$P$ amounts to finding a  factorization of the slack matrix~$\Phi=TS$ into nonnegative matrices~$T$ with~$f$ columns and~$S$ with~$f$ rows. 
In particular, we have derived the following result that essentially is due to Yannakakis~\cite{Yan91} (see also~\cite{FKPT11}). Here, the \emph{nonnegative rank} of a matrix is the minumum number~$f$ such that the matrix can be written as a product of two nonnegative matrices, where the first one has~$f$ columns and the second one has~$f$ rows.
\begin{theorem}
	\label{thm:extComplNonnegRk}
The extension complexity of a  polytope~$P$ is equal to the nonnegative rank of its slack matrix (w.r.t. any  set~$X$ and binding system $Ax\le b$ with $P=\conv{X}=\setDef{x\in\aff{P}}{Ax\le b}$). 	
\end{theorem}
Clearly, the nonnegative rank of a matrix is bounded from below by its usual rank as known from Linear Algebra. There is also quite some interest in the \emph{nonnegative} rank of (not necessarily slack) matrices in general (see, e.g., Cohen and Rothblum~\cite{CR93}).

\section{Fundamental Limits}
\label{sec:fundLims}

\subsection{General Lower Bounds}
Every extension~$Q$ of a polytope~$P$ has at least as many faces as~$P$, as the face lattice of~$P$ can be embedded into the face lattice of~$Q$ (see Sect.~\ref{subsec:facLat}). Since each face is the intersection of some facets, one finds that the extension complexity of a polyhedron with~$\beta$ faces is at least $\log \beta$ (the binary logarithm of~$\beta$). This observation has first been made by Goemans~\cite{Goe09} in order to argue that the extension complexity of the permutahedron~$\perm{n}$ is at least $\Omega(n\log n)$.

Suppose that $\Phi=TS$ is a factorization of a slack matrix~$\Phi$ of the polytope~$P$ into nonnegative matrices~$T$ and~$S$ with columns $t^1,\dots,t^f$ and rows $s^1,\dots,s^f$, respectively. Then we can write $\Phi=\sum_{i=1}^f t^is^i$ as the sum of~$f$ nonnegative matrices of rank one. Calling the set of all non-zero positions of a matrix its \emph{support}, we thus find that the nonnegative factorization $\Phi=TS$ provides a way to  cover the support of~$\Phi$ by~$f$ \emph{rectangles}, i.e., sets of the form $I\times J$, where~$I$ and~$J$ are subsets of the row- and column-indices of~$\Phi$, respectively. Hence, due to Theorem~\ref{thm:extComplNonnegRk}, the minimum number of rectangles by which one can cover the support of~$\Phi$ yields a lower bound (the \emph{rectangle covering bound}) on the extension complexity of~$P$ (Yannakakis~\cite{Yan91}). Actually, the rectangle covering bound dominates the bound discussed in the previous paragraph~\cite{FKPT11}. 
As Yannakakis~\cite{Yan91}  observed furthermore, the logarithm of the rectangle covering bound of a polytope~$P$ is equal to the \emph{nondeterministic communication complexity} (see, e.g., the book of Kushilevitz and Nisan~\cite{KN97}) of the predicate on the pairs $(v,f)$ of vertices~$v$ and facets~$f$ of~$P$ that is true if and only if~$v\not\in f$. 

One can equivalently describe the rectangle covering bound as the minimum number of complete bipartite subgraphs needed to cover the \emph{vertex-facet-non-incidence graph} of the polytope~$P$. A \emph{fooling set} is a subset~$F$ of the edges of this graph such that no two of the edges in~$F$ are contained in a complete bipartite subgraph. Thus every fooling set~$F$ proves that the rectangle covering bound, and hence, the extension complexity of~$P$, is at least~$|F|$.  
For instance, for the $n$-dimensional cube it is not too difficult to come up with a fooling set of size~$2n$, proving that for a cube one cannot do better by allowing extended formulations for the representation. For more details on bounds of this type we refer to~\cite{FKPT11}. 

Unfortunately, all in all the currently known techniques for deriving lower bounds on extension complexities are rather limited and yield mostly quite unsatisfying bounds.


\subsection{The Role of Symmetry}
\label{subsec:symm}
Asking, for instance,  about the extension  complexity of the matching polytope $\allMatchPoly{n}$ defined in the beginning, one finds that not much is known. It might be anything between quadratic and exponential in~$n$. However, in the main part of his striking paper~\cite{Yan91}, Yannakakis established an exponential lower bound on the sizes of \emph{symmetric} extended formulations of~$\allMatchPoly{n}$.  Here, \emph{symmetric} means that the extension polyhedron remains unchanged when renumbering the nodes of the complete graph, or more formally that,
for each permutation~$\pi$ of the edges of the complete graph that is induced by a permutation of its nodes, there is a permutation~$\kappa_{\pi}$ of the variables of the extended formulation that maps the extension polyhedron to itself such that, for every vector~$y$ in the extended space, applying~$\pi$ to the projection of~$y$ yields the same vector as projecting the vector obtained from~$y$ by applying~$\kappa_{\pi}$. 
Indeed, many extended formulations are symmetric in a similar way, for instance the extended formulation of the permutahedron by the Birkhoff-polytope mentioned in the Introduction as well as the extended formulation for the spanning tree polytope discussed in Section~\ref{subsec:spt}.

In order to  state Yannakakis' result more precisely, denote by~$\match{\ell}{n}$ the set of all matchings of cardinality~$\ell$ in the complete graph with~$n$ nodes, and by 
 $\matchPoly{\ell}{n}=\convOp\setDef{\charVec{M}}{M\in\match{\ell}{n}}$ the associated polytope. In particular,  $\matchPoly{n/2}{n}$ is the \emph{perfect-matching-polytope} (for even~$n$).
\begin{theorem}[Yannakakis~\cite{Yan91}]
	\label{thm:YanLBSym}
	For even~$n$, the size of every \emph{symmetric} extended formulation of~$\matchPoly{n/2}{n}$ is at least 
	 $\Omega(\binom{n}{\lfloor (n-2)/4\rfloor})$.
\end{theorem}
Since $\matchPoly{\lfloor n/2\rfloor}{n}$ is (isomorphic to) a face of~$\allMatchPoly{n}$,  one easily derives  the above mentioned exponential lower bound on the sizes of symmetric extended formulations for $\allMatchPoly{n}$ from Theorem~\ref{thm:YanLBSym}. 

At the core of his beautiful  proof of Theorem~\ref{thm:YanLBSym}, Yannakakis shows that, for even~$n$, there is no symmetric extended formulation in equation form (i.e., with equations and nonnegativity constraints only) of~$\matchPoly{n/2}{n}$ of size at most~$\binom{n}{k}$ with $k=\lfloor (n-2)/4\rfloor$. From such a hypothetical extended formulation~$\extForm_1$, he first constructs an extended formulation~$\extForm_2$ in equation form on variables $y_A$ for all matchings~$A$ with $|A|\le k$ such that 
the 0/1-vector valued map~$s^{\star}$ on the vertices of~$\matchPoly{n/2}{n}$ defined by  $s^{\star}(\charVec{M})_A=1$ if and only if $A\subseteq M$ is a \emph{section}  of~$\extForm_2$, i.e., $s^{\star}(x)$ maps every vertex~$x$ to a preimage under the projection of~$\extForm_2$ that is contained in the extension polyhedron.
Then it turns out that an extended formulation like~$\extForm_2$ cannot exist. In fact, for an arbitrary partitioning of the node set into two parts~$V_1$ and~$V_2$ with $|V_1|=2k+1$, one can construct a nonnegative point~$y^{\star}$ in the affine hull of the image of~$s^{\star}$ (thus $y^{\star}$ is contained in the extension polyhedron of~$\extForm_2$ that is defined by equations and nonnegativity constraints only) with $y^{\star}_{\{e\}}=0$ for all edges~$e$ connecting~$V_1$ and~$V_2$, which implies that the projection of the point~$y^{\star}$ violates the inequality $x(\delta(V_1))\ge 1$ that is valid for~$\matchPoly{n/2}{n}$ (since $|V_1|=2k+1$ is odd). The crucial ingredient for constructing~$\extForm_2$ from~$\extForm_1$ is a theorem of Bocherts'~\cite{Boc89} stating that every subgroup~$G$ of  permutations of~$m$ elements that is primitive with $|G|>m!/\lfloor (m+1)/2\rfloor !$ contains all even permutations. Yannakakis constructs a section~$s$ for~$\extForm_1$ for that he can show---by exploiting Bochert's theorem---that there is a nonnegative matrix~$C$ with $s(\charVec{M})=C\cdot s^{\star}(\charVec{M})$ for all $M\in\match{n/2}{n}$, which makes it rather straight forward to construct~$\extForm_2$ from~$\extForm_1$.
 
With respect to the fact that his proof yields an exponential lower bound only for \emph{symmetric} extended formulations, Yannakakis~\cite{Yan91} remarked   ``we do not think that asymmetry helps much'' in constructing small extended formulations of the (perfect) matching polytopes and  stated as an open problem to ``prove that the matching (\ldots) polytopes cannot be expressed by polynomial size LP’s without the symmetry assumption''. As indicated above, today we still do not know whether this is possible. However, at least it turned out recently that requiring symmetry can make a big difference for the smallest possible size of an extended formulation. 

\begin{theorem}[\cite{KPT10}]
	\label{thm:KPT}
	All symmetric extended formulations of~$\matchPoly{\lfloor \log n\rfloor}{n}$ have size at least $n^{\Omega(\log n)}$, while there are polynomial size non-symmetric extended formulations for~$\matchPoly{\lfloor \log n\rfloor}{n}$ (i.e., 
	the extension complexity of~$\matchPoly{\lfloor \log n\rfloor}{n}$ is boun\-ded from above by a polynomial in~$n$).
\end{theorem}

Thus, at least when considering  matchings of size $\lfloor\log n\rfloor$ instead of perfect (or arbitrary) matchings, asymmetry indeed helps much.  

While the proof of the lower bound on the sizes of symmetric extended formulations stated in Theorem~\ref{thm:KPT} is  a  modification of Yannakakis' proof indicated above, the construction of the polynomial size non-symmetric extended formulation of  $\matchPoly{\lfloor \log n\rfloor}{n}$ relies on the principle of disjunctive programming (see Section~\ref{subsec:disPrg}). For an arbitrary coloring~$\zeta$ of the~$n$ nodes of the complete graph with~$2k$ colors, we call a matching~$M$ (with $|M|=k$) \emph{$\zeta$-colorful} if, in each of the~$2k$ color classes, there is exactly one node that is an end-node of  one of the edges from~$M$. Let us denote by~$P_{\zeta}$ the convex hull of the characteristic vectors of $\zeta$-colorful matchings.  The crucial observation is that~$P_{\zeta}$ can be described by $\bigO{2^k+n^2}$ inequalities (as opposed to~$\Omega(2^n)$ inequalities needed to describe the polytope associated with all matchings, see the Introduction). On the other hand, according to a theorem due to  Alon, Yuster, and Zwick~\cite{AYZ95}, there is a family of~$q$ such colorings $\zeta_1,\dots,\zeta_q$ with $q=2^{\bigO{k}}\log n$ such that, for every $2k$-element subset~$W$ of the~$n$ nodes, in at least one of the colorings the nodes from~$W$ receive pairwise different colors. Thus we have $\matchPoly{k}{n}=\conv{P_{\zeta_1}\cup\cdots\cup P_{\zeta_q}}$, and hence (as described in Section~\ref{subsec:disPrg}) we obtain an extended formulation of~$\matchPoly{k}{n}$ of size $2^{\bigO{k}}n^2\log n$, which, for $k=\lfloor\log n\rfloor$, yields the upper bound in Theorem~\ref{thm:KPT}.

Yannakakis~\cite{Yan91} moreover deduced from Theorem~\ref{thm:YanLBSym} that there are  no polynomial size symmetric extended formulations for the traveling salesman polytope (the convex hull of the characteristic vectors of all cycles of lengths~$n$ in the complete graph with~$n$ nodes). Similarly to Theorem~\ref{thm:KPT}, one can also prove that there are no polynomial size symmetric extended formulations for the polytopes associated with cycles of length~$\lfloor\log n\rfloor$, while these polytopes nevertheless have polynomially bounded extension complexity~\cite{KPT10}.

Pashkovich~\cite{Pas11} further extended Yannakakis' techniques in order to prove that every symmetric extended formulation of the  permutahedron~$\perm{n}$ has size at least~$\Omega(n^2)$, showing that the Birkhoff-polytope essentially provides an optimal \emph{symmetric} extension for the permutahedron.

\section{Conclusions}

Many polytopes associated with combinatorial optimization problems can be represented in small, simple, and nice ways as projections of higher dimensional polyhedra. Moreover, though we have not touched this topic here, sometimes such extended formulations are also very helpful in deriving descriptions in the original spaces. What we currently lack are on the one hand more techniques to construct  extended formulations and on the other hand a good understanding of the fundamental limits of such representations. For instance, does every polynomially solvable combinatorial optimization problem admit an extended formulation of polynomial size? We even do not know this for the matching problem. How about the stable set problem in perfect graphs? The best upper bound on the extension complexity of these polytopes for graphs with~$n$ nodes still is~$n^{\bigO{\log n}}$ (Yannakakis~\cite{Yan91}). 

Progress on such questions will eventually shed more light onto the principle possiblities to express combinatorial problems by means of linear constraints. Moreover, the search for extended formulations  yields new modelling ideas some of which may prove to be useful also in practical contexts. In any case,  work on extended formulations can lead into fascinating mathematics. 

 \paragraph{Acknowldgements}
 We are grateful to Sam Burer, Samuel Fiorini, Kanstantsin Pashkovich, Britta Peis, Laurence Wolsey, and Mihalis Yannakakis for  comments on a draft of  this article and to Matthias Walter for 
producing 
Figure~\ref{fig:embedding}.

\end{document}